\documentclass{article}
\usepackage{a4}
\usepackage{amssymb}
\usepackage{amsthm} 
\usepackage{amsmath}
\usepackage{graphics} 
\usepackage{eepic}    
\usepackage[all,arc,poly]{xy}

\newcommand{\nc}{\newcommand}
\nc{\look}{\marginpar{$\bullet$}}

\nc{\Section}{\section}
\nc{\SubSection}{\subsection}

\newtheorem{theo}{Theorem}[section]   
\newtheorem{ddef}[theo]{Definition}
\newtheorem{llem}[theo]{Lemma} 
\newtheorem{oobs}[theo]{Observation} 
\newtheorem{rrem}[theo]{Remark} 
\newtheorem{prop}[theo]{Proposition} 
\newtheorem{ccor}[theo]{Corollary}  
\newtheorem{claim}[theo]{Claim} 
\newtheorem{qquest}[theo]{Question} 
\newtheorem{fact}[theo]{Fact} 
\newtheorem{pprov}[theo]{Proviso}
\newtheorem{eexam}[theo]{Example} 

\nc{\bT}{\begin{theo}} 
\nc{\eT}{\end{theo}}
\nc{\bD}{\begin{ddef} \rm }
\nc{\eD}{\end{ddef}}
\nc{\bC}{\begin{ccor}}
\nc{\eC}{\end{ccor}}
\nc{\bCl}{\begin{claim}}
\nc{\eCl}{\end{claim}}
\nc{\bQ}{\begin{qquest}}
\nc{\eQ}{\end{qquest}}
\nc{\bL}{\begin{llem}}
\nc{\eL}{\end{llem}}
\nc{\bP}{\begin{prop}}
\nc{\eP}{\end{prop}}
\nc{\bR}{\begin{rrem}}
\nc{\eR}{\end{rrem}}
\nc{\bO}{\begin{oobs}}
\nc{\eO}{\end{oobs}}
\nc{\bF}{\begin{fact}}
\nc{\eF}{\end{fact}}
\nc{\bProv}{\begin{pprov}}
\nc{\eProv}{\end{pprov}}
\nc{\bE}{\begin{eexam} \rm }
\nc{\eE}{\end{eexam}}

\nc{\prf}{\begin{proof}}
\nc{\eprf}{\end{proof}}

\renewcommand{\phi}{\varphi}
\renewcommand{\geq}{\geqslant}
\renewcommand{\leq}{\leqslant}

\renewcommand{\subset}{\subseteq}

\nc{\simg}{\sim_{\mathsf{g}}}
\nc{\ssimg}{\approx_{\mathsf{g}}}
\nc{\simC}{\sim_{\ssc \mathsf{\#}}}
\nc{\simCC}[2]{\simC^{#1,#2}}
\nc{\crk}{\mathrm{rk}_{\ssc \mathsf{\#}}}
\nc{\nd}{\mathrm{nd}}
\nc{\ssimC}{\approx_{\mathsf{c}}}
\nc{\ssimCC}{\approx_{\mathsf{cc}}}
\nc{\FO}{\mathsf{FO}}
\nc{\SO}{\mathsf{SO}}
\nc{\GF}{\mathsf{GF}}
\nc{\ML}{\mathsf{ML}}
\nc{\CML}{\mathsf{CML}}
\nc{\gd}{\mathsf{gd}}
\nc{\hg}{\mathsf{hg}}
\nc{\vg}{\mathsf{vg}}
\nc{\free}{\mathrm{free}}
\nc{\dom}{\mathrm{dom}}
\nc{\nn}{\mathsf{non}}
\nc{\NE}{\mathsf{NE}}
\nc{\NF}{\mathsf{NF}}

\newenvironment{romanenumerate}%
{\begin{list}{(\roman{enumi})}{\usecounter{enumi}
\setlength{\labelwidth}{2cm}
\setlength{\itemindent}{0pt}
\setlength{\itemsep}{0.5\itemsep}
\setlength{\topsep}{\itemsep}
\setlength{\parsep}{0pt}
}}{\end{list}}
\nc{\bre}{\begin{romanenumerate}}
\nc{\ere}{\end{romanenumerate}}
\newenvironment{alphaenumerate}%
{\begin{list}{(\alph{enumii})}{\usecounter{enumii}
\setlength{\labelwidth}{2cm}
\setlength{\itemindent}{0pt}
\setlength{\itemsep}{0.5\itemsep}
\setlength{\topsep}{\itemsep}
\setlength{\parsep}{0pt}
}}{\end{list}}
\nc{\bae}{\begin{alphaenumerate}}
\nc{\eae}{\end{alphaenumerate}}
\newenvironment{numenumerate}%
{\begin{list}{(\arabic{enumiii})}{\usecounter{enumiii}
\setlength{\labelwidth}{2cm}
\setlength{\itemindent}{0pt}
\setlength{\itemsep}{0.5\itemsep}
\setlength{\topsep}{\itemsep}
\setlength{\parsep}{0pt}
}}{\end{list}}
\nc{\bne}{\begin{numenumerate}}
\nc{\ene}{\end{numenumerate}}

\nc{\ins}[1]{\bigskip\noindent
\framebox{\begin{minipage}{.95\textwidth} \sloppy \noindent \em #1 \end{minipage}}\bigskip}

\nc{\str}[1]{{\mathfrak{#1}}}
\nc{\brck}[1]{[\![ #1 ]\!]}
\nc{\diam}[1]{\Diamond^{\geq #1}}
\nc{\restr}{\!\restriction\!}
\nc{\G}{\mathbb{G}}
\nc{\HH}{\mathbb{H}}
\nc{\VV}{\mathbb{V}}
\nc{\abar}{\mathbf{a}}
\nc{\bbar}{\mathbf{b}}
\nc{\cbar}{\mathbf{c}}
\nc{\xbar}{\mathbf{x}}
\nc{\ybar}{\mathbf{y}}
\nc{\zbar}{\mathbf{z}}
\nc{\ubar}{\mathbf{u}}
\nc{\vbar}{\mathbf{v}}
\nc{\wbar}{\mathbf{w}}

\nc{\Xbar}{\mathbf{X}}
\nc{\Ybar}{\mathbf{Y}}
\nc{\Zbar}{\mathbf{Z}}
\nc{\Pbar}{\mathbf{P}}
\nc{\nubar}{\mbox{\boldmath $\nu$}}

\nc{\barr}{\begin{array}}
\nc{\earr}{\end{array}}
\nc{\btab}{\begin{tabular}}
\nc{\etab}{\end{tabular}}

\nc{\nothing}{\rule{0em}{1ex}}
\nc{\highnothing}{\rule{0em}{3ex}}
\nc{\hnt}{\highnothing}
\nc{\nt}{\nothing}
\nc{\nnt}{\rule{.1pt}{0pt}}

\nc{\ssc}{\scriptscriptstyle}
\nc{\N}{{\mathbb N}}
\nc{\Z}{{\mathbb Z}}
\nc{\PI}{{\bf I}}
\nc{\PII}{{\bf II}}

\renewcommand{\epsilon}{\varepsilon}

\begin{document}

\title{Graded modal logic and counting bisimulation}

\date{Martin Otto, 2023}

\maketitle

\noindent
\paragraph*{Abstract:}
This note sketches the extension of
the basic characterisation theorems
as the bisimulation-invariant fragment of
first-order logic to modal logic with 
graded modalities and matching adaptation of bisimulation. The result
is closely related to unpublished work in the diploma thesis of 
Rebecca Lukas~\cite{Lukas}, the presentation newly adapted.
We focus on showing \emph{expressive completeness} of graded multi-modal logic
for those first-order properties of pointed Kripke structures 
that are preserved under counting bisimulation equivalence among all
or among just all finite pointed Kripke structures,
similar to the treatment in~\cite{OttoNote,GorankoOtto}. 
The classical version of the characterisation was treated
in~\cite{deRijke}.

\section{Modal back-and-forth with costed counting}
We write $\ML$ for basic modal logic, typically considering its
multi-modal format with finitely many modalities indexed as
$(\Diamond_i)$ with duals $(\Box_i)$, where we think of 
the indices $i$ as ranging over a finite set if \emph{agents}, 
and finitely many basic propositions indexed as $(p_j)$. 
In the mono-modal setting with just one agent  we write 
$\Diamond$ and $\Box$ for the associated modalities. 
Intended models are 
Kripke structures $\str{M} = (W, (E_i), (P_j))$ on sets $W \not=
\emptyset$ of \emph{possible worlds}, with 
accessibility relations $E_i = E_i^\str{M} \subset W\times W$ 
as interpretations of binary relation symbols $E_i$, one for each agent~$i$,
and unary predicates $P_j = P_j^\str{M} \subset W$ as interpretations of unary
relation symbols $P_j$, one for each~$j$.%
\footnote{We drop superscripts to reference $\str{M}$ whenever this 
seems uncritical.}
Formulae are evaluated in \emph{pointed} Kripke structures $\str{M},w$ with a
distinguished world $w \in W$. For any world $u \in W$ we write 
$E_i[u]$ for the set of its immediate $E_i$-successors
\[
E_i[u] = \{ v \in W \colon (u,v) \in E_i \},
\]
which supports the usual interpretation of $\Diamond_i$
according to the inductive clause
\[
\str{M},w \models \Diamond_i \phi \; \mbox{ if } \;
\{ u \in E_i[w] \colon \str{M},u \models \phi \} \not= \emptyset. 
\]

The counting extension of $\ML$, which we denote by
$\CML$, has \emph{graded modalities} 
$\diam{k}_i$ for every $k \in \N \setminus \{ 0 \}$ instead
of just $\Diamond_i$, whose semantics captures that there are 
\emph{at least~$k$ successors} w.r.t.\ accessibility relation $E_i$
 such that \ldots

The defining clause for the semantics of $\diam{k}$ is that 
\[
\str{M},w \models \diam{k} \phi \; \mbox{ if } \;
|\{ u \in E_i[w] \colon \str{M},u \models \phi \}| \geq k.
\]

Using $\brck{\phi}^{\str{M}}$ to denote the extension of 
(the property defined by) $\phi$ in $\str{M}$,  
$\brck{\phi}^{\str{M}} := \{ u \in W \colon \str{M},u \models \phi \}$, 
the above clause can be rewritten as
\[
\brck{\diam{k} \phi }^\str{M}= 
\{ u \in W \colon 
|E_i[u] \cap \brck{\phi}^\str{M} | \geq k
\}.
\] 

In particular, $\diam{1}_i$ is just $\Diamond_i$ so that 
obviously $\ML \subset \CML$ in expressive power. 
On the other hand, the usual standard translation 
of $\ML$ into first-order logic $\FO$ to establish $\ML \subset \FO$
easily extends to $\CML$ (at the expense of more first-order variables 
for the parametrisation of that many distinct $E_i$-successors, 
or alternatively with the use of first-order counting quantifiers
$\exists^{\geq k}$ in a two-variable first-order setting), so that 
\[
\ML \subset \CML \subset \FO,
\]
and these inclusions are easily seen to be strict. For
$\FO$-formulae we denote the usual quantifier rank as in
$\mathrm{qr}(\phi)$. The fragments
$\FO_q := \{ \phi \in \FO \colon \mathrm{qr}(\phi) \leq q \}$
and the induced approximations $\equiv^{\FO_q}$ to elementay
equivalence are defined as usual.

Just as the semantics of basic modal logic $\ML$ is 
invariant under bisimulation equivalence, so 
graded modal logic $\CML$ 
is invariant under the natural refinement 
to \emph{graded bisimulation equivalence}. We give a static 
definition of graded bisimulation \emph{relations} before discussing
the more dynamic intuition based on the associated back-and-forth
(b\&f) game 
underlying this natural notion of structural equivalence. 

\bD
\label{gradbisimdef}
A non-empty binary relation $Z \subset W \times W'$ between  
the universes of two Kripke structures $\str{M} = (W,(E_i),(P_j))$
and $\str{M} = (W',(E_i'),(P_j'))$ is a \emph{graded bisimulation
  relation} 
if it satisfies the following conditions:
\bre
\item[--] \emph{atom equivalence}:
for all $(u,u') \in Z$ and $j$, $u \in P_j \Leftrightarrow u' \in P_j'$.
\item[--] \emph{graded forth}:
for all $(u,u') \in Z$ and $i$, for all $k \geq 1$: 
for pairwise distinct $v_1,\ldots, v_k \in E_i[u]$ there are 
pairwise distinct 
$v_1',\ldots, v_k' \in E_i'[u']$ with $(v_1,v_1'),\ldots, (v_k,v_k') \in Z$.
\item[--] \emph{graded back}:
for all $(u,u') \in Z$ and $i$, for all $k \geq 1$: 
for pairwise distinct $v_1',\ldots, v_k' \in E_i'[u']$
there are pairwise distinct 
$v_1,\ldots, v_k \in E_i[u]$ with $(v_1,v_1'),\ldots, (v_k,v_k') \in Z$.
\ere
We write $\str{M},w \simC \str{M}',w'$, and say that 
these pointed Kripke structures are \emph{graded bisimilar} if 
there is such a graded bisimulation relation with $(w,w') \in Z$.
\eD

The associated graded bisimulation game is played between 
two players, \PI\ and \PII, over the two structures
$\str{M},\str{M}'$. Positions are pairings $(u,u') \in W \times W'$ 
and a single round played from this position allows \PI\ 
to challenge \PII\ in the following exchange of moves
\bre
\item[--] 
\PI\ chooses, for some $i$, a \emph{finite} non-empty subset of
one of $E_i[u]$ or $E_i'[u']$;  
\\
\PII\ must respond with a matching subset 
of $E_i'[u']$ or $E_i[u]$ on the opposite side
of the \emph{same finite size};
\item[--]
\PI\ picks a world in the set proposed by \PII;
\\
\PII\ must respond by a choice of a matching world 
in the set proposed by \PI.
\ere
The pairing between the worlds chosen during the second
stage of this round is the new position in the game.
 
Either player loses in this round if stuck (for a choice or
response), and the second player also loses as soon 
as the position violates atom equivalence.

It is not hard to show that a graded bisimulation relation $Z$
as in Definition~\ref{gradbisimdef} \emph{is} an extensional
formalisation of a 
\emph{winning strategy} in the unbounded (i.e.\ infinite) graded
bisimulation game, which is good for any position $(u,u') \in Z$. 
Atom equivalence for $Z$ ensures
that \PII\ cannot have lost in a position covered by $Z$; and the 
b\&f conditions for $Z$ make sure that \PII\
can respond to any challenge by \PI\ without leaving cover by $Z$.
The winning strategy embodied in $Z$ is in general non-deterministic 
in allowing \PII\ several choices in accordance with $Z$.

Game views of the concept of bisimulations typically suggest the adequate 
\emph{finite approximations} to the a priori essentially infinitary closure
condition implicit in the definition of bisimulation relations: the
natural notion of $\ell$-bisimulation captures the existence of a
winning strategy for \PII\ in the \emph{$\ell$-round game}. It
typically provides exactly the right levels of b\&f equivalence
to match finite levels of logical equivalence based on a gradation by quantifier
rank or nesting depth of modal operators. This approach generates 
all kind of natural variants of the classical
Ehrenfeucht--Fra\"\i ss\'e connection, which produce precise matches
between corresponding levels of b\&f structural equivalences and
logical equivalences.%
\footnote{Recall how the classical Ehrenfeucht--Fra\"\i ss\'e matches
$m$-isomorphy, i.e.\ $m$-round pebble game b\&f equivalence, with
$m$-elemenatry equivalence.}
For basic modal logic $\ML$ as well as for its graded variant $\CML$,
the natural notion of $\ell$-bisimulation $\sim^\ell$ or $\simC^\ell$ 
lies in a limitation to $\ell$ rounds, which for these
logics also corresponds to a limitation of the structural exploration 
to depth $\ell$ from the initial worlds -- which on the logics' side 
in turns matches the natural notions of nesting depths 
for modal operators. This modal nesting depth 
$\nd(\phi)$ is defined for $\CML$ just as for $\ML$, with the crucial
clause in the inductive definition that $\nd(\diam{k}\phi) = \nd(\phi)+1$.

\bD[graded $\ell$-bisimilarity]
\mbox{}\\
Pointed Kripke structures $\str{M},w$ and $\str{M}',w'$
are \emph{graded $\ell$-bisimilar} if player \PII\ has a winning
strategy in the $\ell$-round graded bisimulation game 
on $\str{M},\str{M}'$ starting from 
$(w,w')$; notation: $\str{M},w \simC^\ell \str{M}',w'$.
\eD

In the following we write $\CML_\ell$ for the fragment of formulae of
modal nesting depth up to $\ell$, and 
$\equiv^{\CML_\ell}$ for the notion of equivalence
between pointed Kripke structures that is induced by
indistinguishability in $\CML_\ell$. 

\bL[graded modal Ehrenfeucht--Fra\"\i ss\'e theorem]
\label{EFCML1}
\mbox{}\\
For any two finitely branching pointed Kripke 
structures $\str{M},w$ and $\str{M}',w'$ over a common finite modal
signature, and for all $\ell \in \N$:
\[
\str{M},w \simC^\ell \str{M}',w'
\quad \Longleftrightarrow
\quad
\str{M},w \equiv^{\CML_\ell} \str{M}',w'
\]
\eL

The proof is as usual, but it is worth noting that the corresponding 
modal Ehrenfeucht--Fra\"\i ss\'e theorem, i.e.\ the variant for $\ML$ 
rather than $\CML$, is good across all pointed Kripke structures,
finite branching or not. (It requires the same restriction to finite
signatures, though.) The source of this difference is that
the cardinality of the set of $E_i$-successors
can be pinpointed precisely by a single $\CML_1$-formula if finite, 
but requires an infinite collection of formulae (in $\CML_1$) to be 
characterised as infinite. 

\bE
The following two rooted tree frames of depth $2$,
in a modal signature with a single modality, 
are indistinguishable in $\CML$, but not $\simC^2$-equivalent:
\bae
\item[--]
$\str{M},w$ has of a root node $w$ with children $(u_i)_{i \geq 1}$, 
\\
where 
each $u_i$ has precisely $i$ many children; 
\item[--]
$\str{M}',w'$ has a root node
$w'$ with children $(u_i')_{i \in \N}$, 
\\
where 
each $u_i'$ for $i \geq 1$ has precisely $i$ many children 
\\
and the extra 
$u_0$ has infinitely many children.
\eae
They are $\simC^1$-equivalent though.
\eE

The reason behind this example is that the existence of a successor that has infinitely
many successors is not expressible in $\CML$ (nor is it expressible in 
$\FO$, which precludes easy remedies).
This problem does not occur in restriction to finitely branching 
structures. While such restrictions and problems related to
definability also occur for basic modal logic $\ML$, and are 
familiar from classical modal theory in connection with
Hennessy--Milner phenomena, they here already strike at the level of 
$\simC^\ell$ not just at the level of full $\simC$.

The underlying distinction that really matters for our concerns
lies with the index of the associated equivalence
relations on the class of all (or even just all finite) pointed Kripke
structures. The index of $\sim^\ell$ and $\equiv^{\ML_\ell}$ is finite
as long as the underlying signature 
is finite, for each $\ell \in \N$.  
Contrast this with $\simC^1$, even for just a single modality and 
without propositions. Here we find that the formulae 
$\diam{k} \top$ are pairwise inequivalent, and the family of pointed Kripke
frames consisting of a root node with precisely $k-1$ immediate $E$-successors  
are pairwise $\simC^1$-inequivalent.

To apply a finer gradation than the 
straightforward gradation by just modal nesting depth, we account
for the \emph{cost of counting}. Obviously this is something to take
into account in relation to $\FO$ where any natural rendering of counting
has a cost in terms of quantification.

\bD
[$c$-graded $\ell$-bisimilarity]
\mbox{}\\
For $c \in \N$ consider the $\ell$-round graded bisimulation game with the additional
constraint that the sets $s/s'$ of successor worlds chosen by \PI, and
consequently the responses $s'/s$ given by \PII, must be of
size $\leq c$. Call this the \emph{$c$-graded $\ell$-round game}. 
$\str{M},w$ and $\str{M}',w'$ are \emph{$c$-graded $\ell$-bisimilar},
$\str{M},w \simCC{c}{\ell} \str{M}',w'$, if \PII\
has a winning strategy in this $c$-graded $\ell$-round game
starting from $(w,w')$.
\eD

To match the constraint syntactically, in terms of the grades~$k$ on 
modalities~$\diam{k}$, we define fragments $\CML_{c,\ell}$. 
Essentially we refine the primary gradation of $\CML$ in terms
of nesting depth~$\ell$, which we have in $\CML = \bigcup_\ell \CML_\ell$.

Let the (counting) \emph{rank},
$\crk(\phi)\in \N$,
of formulae $\phi \in \CML$ be defined by induction 
according to
\bre
\item[--] $\crk(\phi) := 0$ for propositional formulae;
\item[--] $\crk(\neg \phi) = \crk(\phi)$;
\item[--] $\crk(\phi_1 \ast \phi_2) = \mathrm{max}(\crk(\phi_1),\crk(\phi_2))$
for $\ast = \wedge,\vee$; 
\item[--] $\crk(\diam{k}_i \phi) := 
\mathrm{max}\bigl(k, \crk(\phi))$.
\ere

The fragment $\CML_{c,\ell}$ comprises all $\CML$-formulae $\phi$
that have counting rank $\crk(\phi) \leq c$ and nesting depth
$\nd(\phi) \leq \ell$. 
Correspondingly we define 
\[
\str{M},w \equiv^{\CML_{c,\ell}}\str{M}',w'
\]
as indistinguishability at level $c,\ell$:
$\str{M},w \models \phi \; \Leftrightarrow\; \str{M}',w' \models \phi$
for all $\phi \in \CML_{c,\ell}$.
Obviously $\CML_{c,\ell} \subset \CML_\ell$
for every $c$, and $\equiv^{\CML_{c,\ell}}$ can only be coarser than 
$\equiv^{\CML_\ell}$
while leading to the same common refinement $\equiv^\CML$
in the limit as
\[
\bigcup_{c,\ell} \CML_{c,\ell} = \CML = \bigcup_\ell \CML_\ell. 
\]

Crucially, for finite modal signatures,
$\CML_{c,\ell}$ is finite up to logical equivalence, for
every combination of $c, \ell \in \N$. 
This is obvious at level $\ell =0$ (by finiteness of the set 
of basic propositions in the signature).
For fixed $c$, the relevant $\diam{k}_i$-constituents at nesting depth
level $\ell+1$ are constrained to finitely many choices for $i$ 
(by finiteness of the set of agents in the signature).
Inductively, the finitely many 
representatives of any constituent subformulae at counting rank level~$c$ 
at nesting depth~$\ell$ can only contribute to finitely many disjunctive normal forms 
involving admissible $\diam{k}_i$-applications at
counting rank level~$c$ and nesting depth $\ell+1$.

\bL[$c$-graded modal Ehrenfeucht--Fra\"\i ss\'e theorem]
\label{EFCML2}
\mbox{}\\
For any two pointed Kripke 
structures $\str{M},w$ and $\str{M}',w'$ over a common finite modal
signature, and for all $\ell \in \N$:
\[
\str{M},w \simCC{c}{\ell}\str{M}',w'
\quad \Longleftrightarrow
\quad
\str{M},w \equiv^{\CML_{c,\ell}}\str{M}',w'
\]
\eL

In comparison with Lemma~\ref{EFCML1}
the restriction to finitely branching structures has been lifted.
We present a detailed proof, which retraces and implicitly reviews
the key idea also for Lemma~\ref{EFCML1} and any other natural variant of the
Ehrenfeucht--Fra\"\i ss\'e connection.

In the following  
we write $[\str{M},w]_{c,\ell}$ for the 
$\simCC{c}{\ell}$-equivalence class of a pointed Kripke structure $\str{M},w$,
\[
[\str{M},w]_{c,\ell} := \{ \str{M}'\!,u'  \colon \; \str{M}',u'  \simCC{c}{\ell} \str{M},u \},
\]
and, relative to a fixed Kripke structure $\str{M}$, just 
$[u]_{c,\ell}$ for the $\simCC{c}{\ell}$-equivalence classes of its worlds $u \in W$:
\[
[u]_{c,\ell} := \{ u' \in W \colon \str{M},u'  \simCC{c}{\ell} \str{M},u \}.
\]

\prf[Proof of Lemma~\ref{EFCML2}.]
To establish the implication from left to right
in the claim of the lemma 
we show that $\crk(\phi) \leq c$, $\nd(\phi) \leq \ell$
and 
$\str{M},w \simCC{c}{\ell}
\str{M}',w'$ imply that 
\[
\str{M},w \models \phi \; \Rightarrow \; 
\str{M}',w' \models \phi. 
\]

This is shown by syntactic induction on $\phi \in \CML$.
The claim is trivial for propositional $\phi$ and obviously compatible
with boolean connectives. For $\phi = \diam{k}_i \psi$, the definition
of counting rank guarantees that $\crk(\phi) \geq k,\crk(\psi)$.

So $\str{M},w \simCC{c}{\ell} \str{M}',w'$ for $\ell \geq \nd(\phi)$ 
and $c \geq \crk(\phi)$ implies 
that \PII\ has, in the first round, a response if \PI\ chooses 
a set $s \subset E_i[w] \cap \brck{\psi}^{\str{M}}$ of size $k$ (such
$s$ exists since $\str{M},w \models \diam{k} \psi$); for this response, say
$s' \subset E_i[w']$ to be adequate, $\PII$ must have, for every 
$u' \in s'$ (that \PI\ could choose in the second part of the first round)
some response $u \in s$ such that $\str{M},u \simCC{c}{\ell -1}
\str{M}',u'$. The inductive hypothesis for $\psi$
implies that $s' \subset \brck{\psi}^{\str{M}'}$ as $\nd(\psi)
\leq \ell-1$ and $\crk(\psi) \leq c$. Thus  
$\str{M}',w' \models \phi$.

For the converse implication, we 
define characteristic formulae $\chi^\ell_{\str{M},w}$ by induction on
$\ell \in \N$ as follows:
\[
\barr{@{}rclr@{}}
\chi^{c,0}_{\str{M},w} &:=& \bigwedge \{p_j \colon w \in P_j^\str{M}  \} 
\wedge \bigwedge \{\neg p_j \colon w \not\in P_j^\str{M}  \} 
\\
\hnt
\chi^{c,\ell+1}_{\str{M},w}&:=& 
\chi^{c,\ell}_{\str{M},w}  
\\
\hnt
&& \!\!\!\wedge\;
\bigwedge \; \{ \;\; \diam{k}_i \chi^{c,\ell}_{\str{M},u} \colon
\; 
u \in E_i[w], 
|E_i[w] \cap [u]_{c,\ell} | \geq k, k \leq c  \} & \mbox{[\emph{forth}]}
\\
\hnt
&& \!\!\!\wedge\;
\bigwedge \; \{ \neg \diam{k}_i \chi^{c,\ell}_{\str{M},u} \colon \; 
u \in E_i[w], |E_i[w] \cap [u]_{c,\ell} | < k \leq c \}
& \mbox{[\emph{back}]}
 \earr
\]

Inductively it is clear that $\chi^{c,\ell}_{-} \in \CML_{c,\ell}$ since
the formally unbounded conjunctions are finite up to logical
equivalence. We now have this, 
without any assumption on branching degree, 
as all $\CML_{c,\ell}$ are finite up to logical equivalence. 
We show that $\chi^{c,\ell}_{\str{M},w}$ defines the
$\simCC{c}{\ell}$-equivalence class 
$[\str{M},w]_{c,\ell}$ of $\str{M},w$.
Obviously $\str{M},w \models \chi^{c,\ell}_{\str{M},w}$; it remains to show that 
\[
\str{M}',w' \models \chi^{c,\ell}_{\str{M},w}
\; \Rightarrow \; 
\str{M},w \simCC{c}{\ell} \str{M}',w'.
\]

This is clear for $\ell = 0$, and we show the claim by inductioin on
$\ell$ for fixed $c$. 
Assuming the claim at level $\ell$, we know that, across all pointed
Kripke structures of the relevant (finite) signature, 
there are finitely many $\simCC{c}{\ell}$-types $t$, each defined by
some characteristic formula $\chi_t^{c,\ell}$. 
To establish the claim at level $\ell+1$ we look at the 
first round in the $c$-graded $(\ell+1)$-round game 
on $\str{M}$ and $\str{M}'$ starting from
$(w,w')$ and give good strategy advice for \PII. We distinguish two
cases, depending on the challenge played by \PI, and invoking
either the \emph{forth} or the \emph{back} constituents in 
$\chi^{c,\ell+1}_{\str{M},w}$.

If \PI's choice is some $s \subset E_i[w]$ of size $|s| \leq c$, 
we partition $s$ into its disjoint constituents
\[
s_t := s \cap \brck{\chi^{c,\ell}_t}^\str{M} 
\]
where $t$ ranges over the finitely many, mutually exclusive, 
$\simCC{c}{\ell}$-types. As $s_t \subset s$ and $k_t := |s_t| \leq c$, 
$\str{M}',w' \models \chi^{c,\ell+1}_{\str{M},w}$ implies in particular 
that, for each $t$,  
\[
\str{M}', w' \models \diam{k_t} \chi_t^{c,\ell}
\]
so that $E_i'[w']$ contains a subset $s'_t$ of size $|s'_t| = k_t =
|s_t|$ of worlds in $\brck{\chi_t^{c,\ell}}^{\str{M}'}$. The
cardinalities of these subsets adds up to $\sum_t |s_t'| = \sum
k_t = |s|$ as they are necessarily disjoint (the
$\chi_t^{c,\ell}$ partition the universe). 
Therefore \PII\ can
respond with $s' := \bigcup_t s'_t \subset E_i'[w']$ to make sure that  
any choice of $u' \in s'$ by \PI\ allows her a matching choice 
of some $u \in s$, viz.\ picking $u \in s_t$ if $u' \in s'_t$. That
guarantees $\str{M},u \simCC{c}{\ell} \str{M}',u'$, and hence a win 
in the remaining game.

If \PI's choice is some $s' \subset E_i'[w']$ of size $|s| \leq c$, 
we analogously partition $s'$ into its disjoint constituents
\[
s_t' := s' \cap \brck{\chi^{c,\ell}_t}^{\str{M}'} 
\]
and let $k_t := |s'_t| \leq c$.
Now $\str{M}',w' \models \chi^{c,\ell+1}_{\str{M},w}$
implies that $\str{M},w \models 
\diam{k_t}_i \chi^{c,\ell}_t$: otherwise the negation of this formula
would be a conjunct in the \emph{back}-part 
of $\chi^{c,\ell+1}_{\str{M},w}$ whereas  $\str{M}',w' \models \diam{k_t}_i \chi^{c,\ell}_t$.
Analogous to the above, the union of size $k_t$-subsets
of $s_t \subset E_i[w] \cap \brck{\chi^{c,\ell}_t}^{\str{M}}$ serves as a
safe response for \PII.
\eprf

\bC
\label{NFcor}
Across all (finite or infinite) pointed Kripke structures 
over any fixed finite modal signature, every 
formula $\phi \in \CML_{c,\ell}$ is logically equivalent 
to a disjunction over 
mutually incompatible characteristic formulae $\chi_t^{c,\ell}$ 
from the finite collection of such characteristic formulae $\chi_t^{c,\ell}$ 
for all $\simCC{c}{\ell}$-types $t$. In particular, 
$\equiv^{\CML_{c,\ell}}$ and $\simCC{c}{\ell}$ have 
finite index for all $c,\ell \in \N$ over the class of 
all pointed Kripke structures of fixed finite modal signature.
\eC

\bC
\label{EFcor2}
The following are equivalent for any class $\mathcal{C}$ of pointed Kripke
structures over a finite modal signature and $c,\ell \in \N$:
\bre
\item
$\mathcal{C}$ is definable by a formula $\phi \in \CML_{c,\ell}$, 
i.e.\
\\
$\mathcal{C} = \mathrm{Mod}(\phi)$ for some $\phi \in \CML$ of
$\crk(\phi) \leq c$ and $\nd(\phi) \leq \ell$;
\item
$\mathcal{C}$ is closed under $\simCC{c}{\ell}$.
\ere
The same is true of any class 
$\mathcal{C}$ of \emph{finite} pointed Kripke
structures over a finite modal signature: 
$\mathcal{C} = \mathrm{FMod}(\phi)$ for some $\phi \in \CML_{c,\ell}$ if,
and only if, $\mathcal{C}$ is closed under $\simCC{c}{\ell}$ within the
class of all finite pointed Kripke structures over the same finite signature.
\eC

So, both in the sense of classical and in the sense of finite model 
theory, definability in $\CML$ is equivalent with $\simCC{c}{\ell}$-invariance 
for some finite level $c,\ell \in \N$.  

Just as for basic modal logic $\ML$ in relation to ordinary
bisimulation $\sim$ and its finite approximations $\sim^\ell$, 
this does not immediately relate definability to \emph{bisimulation invariance}. 
The missing link in both cases is that invariance under the proper 
\emph{infinitary} bisimulation equivalence ($\sim$ or $\simC$) 
is a much weaker (!) condition than invariance under any one of its 
much coarser finite approximations 
($\sim^\ell$ for $\sim$, or either gradation $\simC^\ell$ or
$\simCC{c}{\ell}$ for $\simC$).

Just as in the case of basic modal logic and $\sim$,
the crux for a characterisation of $\CML$ as the 
$\simC$-invariant fragment of $\FO$ therefore lies in establishing 
that for $\phi \in \FO$, 
$\simC$-invariance implies $\simCC{c}{\ell}$-invariance for some 
$c,\ell \in \N$. The crucial difference in these matters arises 
from the fact that the corresponding assertion for $\simC^\ell$
would not suffice: there is no obvious analogue for Corollary~\ref{EFcor2}
with $\simCC{c}{\ell}$ replaced by $\simC^\ell$.

\section{Characterisation theorems through upgrading}

We separately state the two versions of the desired characterisation theorem,
one for the setting of classical model theory and one for finite model theory.

\bT
\label{chathm1}
The following are equivalent for any $\phi(x) \in \FO$
in (the first-order counterpart of) a modal signature:
\bre
\item
$\phi$ is invariant under $\simC$ within the class of all
pointed Kripke structures over the modal signature of $\phi$;
\item
$\phi$ is expressible in graded modal logic $\CML$: 
$\phi \equiv \phi'$ for some $\phi \in \CML$. 
\ere 
\eT

\bT
\label{chathm2}
The following are equivalent for any $\phi(x) \in \FO$
in (the first-order counterpart of) a modal signature:
\bre
\item
$\phi$ is invariant under $\simC$ within the class of all finite 
pointed Kripke structures over the modal signature of $\phi$;
\item
over the class of all finite 
pointed Kripke structures $\phi$ is expressible in 
graded modal logic $\CML$: 
$\phi \equiv_{\mathrm{fin}} \phi'$ for some $\phi \in \CML$. 
\ere 
\eT

In both cases, the implications from~(ii) to~(i) are direct
consequences of the Ehrenfeucht--Fra\"\i ss\'e analysis above.
In both cases the essential part is the \emph{expressive completeness} 
claim of $\CML$ for $\simC$-invariant first-order properties. And in
both cases this claim is equivalent with the assertion that for
first-order definable properties of (finite) pointed Kripke structures 
over a finite modal signature 
\[
(\ddagger) \qquad
\mbox{$\simC$-invariance}
\quad \mbox{ implies } \quad
\mbox{ $\simCC{c}{\ell}$-invariance for some $c,\ell$.}
\]

As I argued elsewhere in connection with the analogue for 
ordinary bisimulation $\sim$ and basic modal logic $\ML$, 
this can be seen as a special \emph{compactness property} 
for $\FO$, and one that unlike full compactness does not fail 
in restriction to finite models.

Just like the very similar technique originally proposed for plain
$\ML$ and $\sim$ in~\cite{OttoNote,GorankoOtto}, and later extended 
to several variants e.g.\ in~\cite{OttoAPAL04,DawarOttoAPAL09}, 
the approach via an upgrading argument for $(\ddagger)$ 
has the advantage of establishing $(\ddagger)$ simultaneously 
for the classical and the finite model theory version.

\medskip
First-order logic is known to satisfy strong locality properties
over relational structures, as expressed in Gaifman's theorem
(cf.\ eg.~\cite{EbbinghausFlu99}).
Recall that a first-order formula $\phi(x)$ (in a single free varaible
$x$) is \emph{$\ell$-local} if its semantics only concerns
the $\ell$-neighbourhood $N^\ell(w)$ of the element $w$ assigned to $x$ 
in a relational structure $\str{M}$. Since we are here dealing with
relational structures in a modal signature, 
the $\ell$-neighbourhood $N^\ell(w)$ of a world $w$ in 
$\str{M} = (W, (E_i^{\str{M}}), (P_j^{\str{M}}))$ 
consists of all worlds $v \in W$ that are at graph distance 
up to $\ell$ from $w$ in the \emph{undirected} 
graph induced by the union of the symmetrisations of all the $E_i^{\str{M}}$.
In this context $\phi(x)$ is $\ell$-local if, for all pointed Kripke
structures $\str{M},w$
\[
\str{M},w \models \phi \; \Leftrightarrow \; 
\str{M} \restr N^\ell(w),w \models \phi, 
\]
where $\str{M} \restr N^\ell(w)$ refers to the substructure 
of $\str{M}$ that is induced on $N^\ell(w) \subset W$.

In preparation for the upgrading argument towards $(\ddagger)$ 
consider the evaluation of any particular 
$\ell$-local first-order formula $\phi^\ell(x)$ 
at some world $w$ in $\str{M}$ for which the induced substructure 
$\str{M}\restr N^\ell(w)$ is \emph{tree-like} 
in the sense that the symmetrisations of the $E_i^{\str{M}}$  
\bre
\item[--] are pairwise disjoint, and
\item[--] their union is acyclic
\ere
in restriction to $N^\ell(w)$. In particular the undirected graph 
induced by the union of the accessibility relations $E_i$ in
restriction to $N^\ell(w)$ is an undirected tree of depth at most
$\ell$ from its root $w$, in the graph-theoretic sense.

We call $\str{M},w$ \emph{rooted tree-like to depth~$\ell$}
if $\str{M}\restr N^\ell(w),w$ is tree-like and the 
$E_i$ themselves happen to be uniformly 
directed away from the root $w$ in $\str{M}\restr N^\ell(w),w$. 

For rooted tree-like $\str{M}\restr N^\ell(w),w$, whether or not
\[
\str{M}\restr N^\ell(w),w \models \phi^\ell(x)
\]
is fully determined by the $\simC^\ell$-type of $\str{M},w$ 
and even by its $\simCC{c}{\ell}$-type provided $c$
is chosen large enough w.r.t.\ the quantifier rank of $q = \mathrm{qr}(\phi)$
and the number of $\equiv^{\FO_q}$-types (cf.\ Lemma~\ref{FOtreelem} below).
The condition that all $E_i$ be directed away from the root
is essential because no degree of $\simC^\ell$-equivalence can 
control in-degrees for $E_i$-edges, while $\FO$ can; and the same goes
for disjointness of the $E_i$ and absence of cycles.

\bL
\label{rootedtreelikelem}
Any $\phi(x) \in\FO_q$ 
that is invariant under $\simC$ on all (or just all finite)
pointed Kripke structures, is 
$\ell$-local for $\ell = 2^q -1$
in restriction to all (or just all finite)
pointed Kripke structures that are rooted tree-like to depth $\ell$.
For such $\str{M},w$, where $\str{M}\restr N^\ell(w),w$ is rooted tree-like:
\[
\str{M},w \models \phi 
\quad \Leftrightarrow 
\quad 
\str{M}\restr N^\ell(w),w \models \phi.
\] 
\eL

The proof is strictly analogous to that given for $\sim$-invariant
$\phi$ in~\cite{OttoNote,GorankoOtto}. Essentially,  
the invariance of $\phi$ under \emph{disjoint unions} with other 
component structures (a fundamental consequence of 
graded as well as plain bisimulation invariance)
reduces the claim to an elementary b\&f 
argument for the $q$-round Ehrenfeucht--Fra\"\i ss\'e game  
on 
\[
\barr{rrclcccl}
& \str{M}',w &:=& 
 q \otimes \str{M}
& \oplus 
& \str{M},w 
& \oplus &  q \otimes \str{M}\restr N^\ell(w)
\\
\hnt
\mbox{versus}
&\str{M}'',w &:=& 
q \otimes \str{M}
& \oplus &\str{M}\restr N^\ell(w) ,w
& \oplus & q \otimes \str{M}\restr N^\ell(w)
\earr
\]
where $\oplus$ stands for disjoint union and 
$q \otimes \str{M}$ for the disjoint union of $q$ many 
copies of $\str{M}$. 

Unlike the treatment for basic $\ML$ cited above, there is
not the immediate shortcut 
that $\ell$-locality together with $\simC$-invariance would imply 
 $\simCC{c}{\ell}$-invariance as required for $(\ddagger)$; and the
obvious analogue for $\simC^\ell$ is not good enough as it does not
feed directly into Corollary~\ref{EFcor2}; 
this is what makes the situation slightly more interesting.
Rather, we need the following simple fact about $\FO$ on tree
structures. The proof is again by a standard b\&f argument about 
the classical $q$-round Ehrenfeucht--Fra\"\i ss\'e game for $\FO$
to establish $\simeq_q$.

\bL
\label{FOtreelem}
For any fixed finite signature and $q,\ell \in \N$ 
there is $c \in \N$ s.t.\
for any two rooted tree-like structures $\str{M}\restr N^\ell(w),w$ 
and $\str{M}' \restr N^\ell(w'),w'$:
\[
\str{M}\restr N^\ell(w),w  \simCC{c}{\ell} \str{M}' \restr
N^\ell(w'), w'
\;\;
\Rightarrow
\;\;
\str{M} \restr N^\ell(w),w 
\equiv^{\FO_q} \str{M}' \restr N^\ell(w'),w',
\]
i.e.\ for structures that are rooted tree-like to depth~$\ell$, the 
$\simCC{c}{\ell}$-type fully determines the 
$\FO_q$-type of their restrictions to depth $\ell$, and hence 
the truth value of any $\ell$-local $\phi(x) \in \FO_q$.
\eL

\bL[upgrading lemma]
\label{upgradelem}
\mbox{}\\
For $\phi(x) \in \FO$ that is invariant under $\simC$ over the class
of all (or just all finite) pointed Kripke structures, there are 
$c,\ell \in \N$ s.t.\ $\phi(x)$ is invariant under
$\simCC{c}{\ell}$ over the class of all (or just all finite) pointed 
Kripke structures. Just as for basic modal logic, the optimal choice
for $\ell$ is $\ell = 2^q -1$ where $q$ is the
quantifier rank of $\phi$.
\eL

This then establishes $(\ddagger)$ and yields the
characterisation theorems.

\prf[Proof of the lemma.]
Let $q := \mathrm{qr}(\phi)$ and $\ell := 2^q-1$. For  
pointed Kripke structures
$\str{M},w$ and $\str{M}',w'$ over the finite signature of $\phi$, 
let $\str{M}^\ast$ and $\str{M}'\nt^\ast$ be the 
partial directed tree unravellings to 
depth~$\ell+1$ from roots $w$ and $w'$, 
merged with copies of the original structures for 
the continuation beyond that depth; this preserves $\simC$ 
and, where relevant, finiteness of the models in question. 
Let $c$ be as in Lemma~\ref{FOtreelem} 
for this situation, and assume now that 
$\str{M},w$ and $\str{M}',w'$, and hence also
$\str{M}^\ast,w$ and $\str{M}'\nt^\ast,w'$, are $\simCC{c}{\ell}$-equivalent 
for these values. 
By Lemma~\ref{rootedtreelikelem}, $\phi(x)$ is $\ell$-local on
$\str{M}^\ast,w$ and $\str{M}'\nt^\ast,w'$, and by
Lemma~\ref{FOtreelem},  
$\str{M}^\ast\restr N^\ell(w),w$ and $\str{M}'\nt^\ast\restr
N^\ell(w'),w'$ agree on $\phi(x)$ since they are $\FO_q$-equivalent.
So
\[
\barr{rcll}
\str{M},w \models \phi(x) &\Leftrightarrow&
\str{M}^\ast,w \models \phi(x) & \mbox{ (by $\simC$-invariance)}
\\
&\Leftrightarrow&
\str{M}^\ast\restr N^\ell(w) ,w \models \phi(x) & \mbox{ (by
  $\ell$-locality, Lemma~\ref{rootedtreelikelem})}
\\
&\Leftrightarrow&
\str{M}'\nt^\ast\restr N^\ell(w') ,w' \models \phi(x) & \mbox{ (by
  choice of $c$,  Lemma~\ref{FOtreelem})}
\\
&\Leftrightarrow&
\str{M}'\nt^\ast ,w' \models \phi(x) & \mbox{ (by
  $\ell$-locality, Lemma~\ref{rootedtreelikelem})}
\\
&\Leftrightarrow&
\str{M}',w' \models \phi(x) & \mbox{ (by $\simC$-invariance)}
\earr
\]
which establishes $\simCC{c}{\ell}$-invariance of $\phi(x)$. 
\eprf

\paragraph*{Remarks.} 
The classical characterisation of $\CML$  
of Theorem~\ref{chathm1} was obtained, in close analogy with
van~Benthem's classical, compactness- and saturation-based 
treatment for plain $\ML$ in~\cite{Benthem83}, by de~Rijke in~\cite{deRijke}.
The finite model theory version in Theorem~\ref{chathm2} is there
stated as an open problem, without reference to the finite model theory
analogue for plain
$\ML$ due to Rosen~\cite{Rosen}. I am not aware of published
work filling that gap, and here took the opportunity to adapt the 
combined and very elementary treatment of  
the van~Benthem--Rosen results from~\cite{OttoNote,GorankoOtto}
to $\CML$ in the most straightforward manner.
In following the well established route through an upgrading
argument that can be made to work in both classical and finite 
model theory context, it also retains much of the
variability that this approach has shown in applications to 
other formats for the underlying modal core logic. In particular, natural
analogous characterisations obtain for some variants involving
global~two-way modalities with counting thresholds, 
and possibly restrictions to special classes of frames, 
similar to some of the corresponding modification treated 
e.g.\ in~\cite{OttoAPAL04,DawarOttoAPAL09}. Some such were 
already treated in Rebecca Lukas' diploma thesis~\cite{Lukas}.

\end{document}